\documentclass[preprint,12pt]{elsarticle}

\usepackage{amssymb,amsmath,lineno,bm,graphicx,multirow,algorithm,algpseudocode,gensymb,textcomp}

\usepackage{hyperref}
\hypersetup{colorlinks=true,
            linkcolor=blue,
            anchorcolor=blue,
            citecolor=blue}

\journal{Journal of Computational Physics}

\begin{document}

\begin{frontmatter}

\title{Reduced-Order Solution for Rarefied Gas Flow by Proper Generalised Decomposition}

\author[label1,label2]{Wei Su\corref{cor1}} 
\author[label3]{Xi Zou}

\affiliation[label1]{organization={Division of Emerging Interdisciplinary Areas, The Hong Kong University of Science and Technology},
            addressline={Clear Water Bay, Kowloon}, 
            city={Hong Kong},
            country={China}}
\affiliation[label2]{organization={Department of Mathematics, The Hong Kong University of Science and Technology},
            addressline={Clear Water Bay, Kowloon}, 
            city={Hong Kong},
            country={China}}

\affiliation[label3]{organization={Zienkiewicz Institute for Modelling, Data and AI, Swansea University},
            city={Swansea},
            postcode={SA1 8EN},
            country={United Kingdom}}

\cortext[cor1]{weisu@ust.hk}

\begin{abstract}
Modelling rarefied gas flow via the Boltzmann equation plays a vital role in many areas. Due to the high dimensionality of this kinetic equation and the coexistence of multiple characteristic scales in the transport processes, conventional solution strategies incur prohibitively high computational costs and are inadequate for rapid response for parametric analysis and optimisation loops in engineering design simulations. This paper proposes an \textit{a priori} reduced-order method based on the proper generalised decomposition to solve the high-dimensional, parametrised Shakhov kinetic model equation. This method reduces the original problem into a few low-dimensional problem by formulating separated representations for the low-rank solution, as well as data and operators in the equation, thereby overcoming the curse of dimensionality. Furthermore, a general solution can be calculated once and for all in the whole range of the rarefaction parameter, enabling fast and multiple queries to a specific solution at any point in the parameter space. Numerical examples are presented to demonstrate the capability of the method to simulate rarefied gas flow with high accuracy and significant reduction in CPU time and memory requirements.
\end{abstract}



\begin{keyword}
reduced-order modelling, proper generalised decomposition, Boltzmann equation, parametrised rarefied gas flow

\end{keyword}

\end{frontmatter}


\section{Introduction}
\label{sec1}

Kinetic theory has demonstrated its practical significance since the 1950s in describing the dynamics of rarefied gas flows encountered in various engineering applications, such as nano-/microelectromechanical systems, high-altitude flights, unconventional natural gas production, vacuum science, etc. The centre of the theory is the Boltzmann equation, which determines the thermofluid properties of a gaseous system by providing evolution information on the probability distribution of gas molecules~\cite{Harris1971}. The Boltzmann equation for monatomic gases reads
\begin{equation*}
    \frac{\partial f}{\partial t}+\bm{v}'\cdot\frac{\partial f}{\partial\bm{x}'}+\bm{F}\cdot\frac{\partial f}{\partial\bm{v}'}=\mathcal{C}\left(f,f\right),
\end{equation*}
where $f\left(t,\bm{x}',\bm{v}'\right)$ is the one-particle velocity distribution function, which is a function of the time $t$, spatial position $\bm{x}'$ and molecular velocity $\bm{v}'$; $\bm{F}$ is external driven force; and $\mathcal{C}$ is the Boltzmann collision integral operator, representing the variation rate of the velocity distribution function due to molecular collisions. Methods for solving the Boltzmann equation are generally categorised as stochastic and deterministic. Stochastic methods use simulation particles to trace the movements and collisions of molecules statistically~\cite{Bird1994}, whereas deterministic methods discretise the velocity distribution function in all independent variables and solve the resultant system algebraically~\cite{Dimarco2014}. Both methods are expensive in computational time and memory requirements, mainly because of the high dimensionality of the equation and the multiscale nature of rarefied gas transport. 

The Boltzmann equation is seven-dimensional. To solve it deterministically, one first discretises the velocity space by $N_{\text{v}}$ discrete points, resulting in $N_{\text{v}}$ partial differential equations (PDE) continuous in time and spatial space that can be solved by finite difference, finite volume, or finite element techniques combing with a time-stepping scheme. Typically, $N_{\text{v}}$ is around $10^4\sim10^6$ for a satisfactory solution, leading to a large number of degrees of freedom (DoF). Meanwhile, a stochastic method must trace many gas molecules in space and time and cope with the variance reduction issue.

In practical engineering circumstances, rarefied gas flows often occur in multiple spatio-temporal scales, characterised by a significant separation of the Knudsen number $\mathtt{Kn}$, which is defined as the ratio of the mean free path or the mean free time of gas molecules to a relevant process length or time scale. Therefore, parametric analysis is necessary to simulate gas flows over a wide range of Knudsen numbers. Due to the complexity of the Boltzmann collision operator, the source iteration scheme is required even when solving stationary and linear problems~\cite{ADAMS20023}. The multiscale nature increases the computational cost, as the convergence of the source iteration is largely slowed when the Knudsen number is small~\cite{SU2020109245}. Furthermore, if a numerical method is not asymptotic preserving, the spatial cell size and time step should be small enough to resolve the molecular mean free path and time to reduce numerical dissipation~\cite{Guo2023UP}. Consequently, inquiring solutions for varying Knudsen numbers is costly because it relies on many expensive, individual simulations.  

Novel and efficient solution strategies are needed for the Boltzmann equation to expedite engineering simulation and design, where fast response is commonly required in parametric analysis and optimisation loops. In recent years, reduced-order modelling has gained increasing interest in reducing the computational complexity of high-dimensional and/or parametrised PDEs, where parameters can be considered as extra coordinates~\cite{Quarteroni2013}. The basic idea is to project a solution into a low-dimensional space spanned by a reduced basis constructed in devised algorithms. The resultant solution is a low-rank or reduced-order solution, in contrast to the full-rank or full-order one. A classic example of model order reduction is digital image compression via singular value decomposition (SVD), which allows for the retention of the most important features of a matrix by truncating the redundant information. To date, popular reduced-order modelling approaches include the reduced basis method~\cite{Rozza2008}, proper orthogonal decomposition~\cite{Berkooz1993}, proper generalised decomposition (PGD)~\cite{PGD2011}, and the dynamic low-rank approximation (DLRA)~\cite{Koch2007}. Among these methods, PGD and DLRA build the reduced basis without \textit{a prior} knowledge of full-rank solutions; thus, they are particularly attractive when it is intractable to compute even one such solution.

The PGD technique features a separated representation: the low-rank solution is approximated as a sum of a few functional products, called \emph{modes}, with each function depending on fewer coordinates of a high-dimensional and/or parametrised problem~\cite{PGDprimer}. As a result, the growth of DoF with respect to the problem's dimension is reduced from an exponential increase to a linear increase at most, breaking the so-called \emph{curse of dimensionality}. The PGD computational process comprises two stages: an \textit{offline} phase calculates modes, once and for all, providing a generalised solution (or the so-called \emph{computational vademecum}) containing all possible solutions for a given domain of all coordinates; and an \textit{online} phase that accesses a specific solution simply via a simple reconstruction from the separated representation, enabling a fast response in parametric analysis. PGD has been used in structural optimisation~\cite{structral}, geometry parametrisation~\cite{SEVILLA2020112631}, inverse parameter identification~\cite{NADAL2015113}, rheology~\cite{CHINESTA2011578}, biomechanics~\cite{Zou2018}, power supply systems~\cite{Power}, fluid dynamics~\cite{TSIOLAKIS2022110802}, and neutron transport~\cite{DOMINESEY2023112137}. The DLRA model finds a low-rank solution by orthogonally projecting the governing equation onto the tangent space of the low-rank solution manifold. Unlike PGD, which constructs the low-rank approximation, DLRA evolves a low-rank solution conventionally (perhaps unnecessarily) over time. DLRA has been applied to the solution of high-dimensional kinetic equations such as the Vlasov equation~\cite{Einkemmer2018}, the radiation transport equation~\cite{PENG2020109735}, and the Boltzmann equation for rarefied gases~\cite{BoltzmannDLRA}.

In this work, the authors aim to explore, for the first time, the PGD method's ability as an efficient solver for the high-dimensional Boltzmann equation and rarefied gas flows parametrised with the rarefaction. Algorithms are developed to transform the solution into a few low-dimensional problems, allowing flow fields to be obtained within one minute regardless of the Knudsen number. For a parametrised problem, the Knudsen number serves as an additional coordinate. The generalised solution of the rarefied gas flow over the whole range of Knudsen numbers is constructed, allowing access to the result at any Knudsen number in real time. In addition, the required memory can also be significantly reduced. The remainder of the paper is arranged as follows. The rarefied gas flow and its kinetic description considered in this work are presented in \S\ref{Sec:Statement}. The separated representation of the low-rank solution and the PGD algorithm to find the solution are detailed in \S\ref{Sec:PGD}. Numerical results are given in \S\ref{Sec:Result}. The manuscript is closed with some conclusions in \S\ref{Sec:conclusion}.

\section{Statement of the Problem}\label{Sec:Statement}

Consider a rarefied gas flowing through a straight channel along $z'$-axis with an arbitrary cross section in the $x'$-$y'$ coordinate plane. The channel connects two reservoirs that contain the same gas and maintain a pressure $P_1$ and a temperature $T_1$ in one and $P_2$ and $T_2$ in the other, with $P_1\leqslant P_2$ and $T_1<T_2$. Due to the pressure and temperature gradients, the gas flow combines a Poiseuille flow from the high to the low pressure reservoir and a thermal creep flow from the cold to the hot reservoir. 

The end effect can be neglected when the channel length is sufficiently larger than its cross-sectional dimension $H$. Furthermore, the dimensionless local pressure and temperature gradients at any cross section along the flow direction $z'$ are always much less than 1~\cite{Sharipov1997}; that is
\begin{equation}
    \xi_P=\frac{H}{P}\frac{\mathrm{d}P}{\mathrm{d}z'}\ll1,\quad \xi_T=\frac{H}{T}\frac{\mathrm{d}T}{\mathrm{d}z'}\ll1,
\end{equation}
where $P=P(z')$ and $T=T(z')$ are the local pressure and temperature, respectively. Due to the small gradients, the velocity distribution function to describe the dynamics of the gas system can be linearised as
\begin{equation}
    f=\frac{n}{v^3_{\text{m}}}\left\{f_{\text{eq}}+\xi_P\left[h_P+zf_{\text{eq}}\right]+\xi_T\left[h_T+z\left(\bm{v}^2-\frac{5}{2}\right)f_{\text{eq}}\right]\right\},
\end{equation}
where $n$ is the gas number density and $v_{\text{m}}=\sqrt{2k_{\text{B}}T/m}$ is the most probable molecular velocity, with $k_{\text{B}}$ being the Boltzmann constant and $m$ the molecular mass. $h_a=h_a\left(x,y,\bm{v}\right)$, $a=\{P,T\}$ are dimensionless perturbed velocity distribution functions responsible for the pressure and temperature gradients, respectively. Note that $\bm{x}=\left(x,y,z\right)$ is the spatial coordinates normalised by $H$, and $\bm{v}=\left(v_x,v_y,v_z\right)$ is the molecular velocity normalised by $v_{\text{m}}$. $f_{\text{eq}}$ is the dimensionless global equilibrium distribution, read
\begin{equation}
    f_{\text{eq}}=\frac{1}{\pi^{3/2}}\exp\left(-\bm{v}^2\right).
\end{equation}

In this work, it is assumed that variations of $h_a$ are governed by the steady-state Shakhov kinetic model equation subjected to the fully diffuse boundary condition, resulting in a first order PDE problem: find $h_a(x,y,\bm{v})$ such that
\begin{equation}\label{Shakhov}
\begin{aligned}
    v_x\frac{\partial h_a}{\partial x}+v_y\frac{\partial h_a}{\partial y}=\delta\left[2u_av_zf_{\text{eq}}+\frac{4}{15}q_av_z\left(\bm{v}^2-\frac{5}{2}\right)f_{\text{eq}}-h_a\right]-s_a,\\ 
    \left(x,y\right)\in\Omega,\\
    h_a=0,\quad \left(x,y\right)\in\partial\Omega,\ \bm{v}\cdot\bm{n}<0,
\end{aligned}
\end{equation}
where
\begin{equation}
    s_a=\begin{cases}
        v_zf_{\text{eq}},&a=P,\\
        v_z\left(\bm{v}^2-\dfrac{5}{2}\right)f_{\text{eq}},\quad&a=T.
    \end{cases}
\end{equation}
In Eq.~\eqref{Shakhov}, $\Omega$ represents the region of the channel cross section and $\partial\Omega$ is its boundary; $\bm{n}$ is the outward unit normal vector at the boundary. $\delta$ is the rarefaction parameter, which is inversely proportional to the Knudsen number
\begin{equation}
    \delta=\frac{PH}{\mu v_{\text{m}}}\propto\frac{1}{\mathtt{Kn}},
\end{equation}
where $\mu$ is the gas viscosity. At the macroscale, the dimensionless flow velocity $u_a$ and the heat flux $q_a$ are obtained directly from integrating the perturbed distribution functions through velocity moments, as
\begin{equation}\label{moments}
    u_a=\int v_zh_a\mathrm{d}\bm{v},\quad q_a=\int v_z\left(\bm{v}^2-\frac{5}{2}\right)h_a\mathrm{d}\bm{v}.
\end{equation}
Dimensionless flow rates are also of interest. The one for the pressure-driven flow is computed as
\begin{equation}
    G_P=-2\int_{\Omega} u_P\mathrm{d}x\mathrm{d}y,
\end{equation}
while that of the temperature-driven flow can be obtained according to the Onsager-Casimir relation~\cite{Sharipov1999}
\begin{equation}
    G_T=2\int_{\Omega} q_P\mathrm{d}x\mathrm{d}y.
\end{equation}

For deterministic solutions of the kinetic equation~\eqref{Shakhov}, the so-called discrete velocity method is commonly used, where the continuous molecular velocity is replaced by $N_{\text{v}}=N_x\times N_y\times N_z$ discrete points $\bm{v}^{j}=\left(v^{j_1}_x,v^{j_2}_y,v^{j_3}_z\right)$ with $j_1=\{1,\dots,N_x\}$, $j_2=\{1,\dots,N_y\}$, and $j_3=\{1,\dots,N_z\}$. Then, the governing equations for the unknown distribution functions $h^{j}_a=h_a\left(x,y,\bm{v}^{j}\right)$ at discrete velocities are solved using a traditional numerical scheme for first order PDEs.

In this work, the nodal discontinuous Galerkin finite element method is used. The discrete distribution functions are approximated in the piecewise polynomial space of complete degree at most $p$, the total DoF to approximate the system is $N_{\text{s}}\times N_{\text{v}}$, where $N_{\text{s}}=N_{\text{el}}\times N_p$ with $N_{\text{el}}$ being the number of finite elements partitioning the computation domain $\Omega$ and $N_p$ the internal DoF in one element~\cite{SU2015123}, for instance, $N_p=\left(p+1\right)\left(p+2\right)/2$ for a triangular element. The flow velocity and heat flux $u_a$ and $q_a$ on the right-hand side of Eq.~\eqref{Shakhov} are evaluated by some quadrature of Eq.~\eqref{moments}, involving all the discrete velocity distribution functions. To avoid an extensive linear system, $h^{j}_a$ is obtained iteratively, where $u_a$ and $q_a$ are estimated from the distribution functions in the previous iteration step and are treated as known when finding the distribution functions in the current step~\cite{SUImplicitDG}. Consequently, let $N^{\text{fr}}_{\text{itr}}$ be the number of iterations for a convergent solution, the overall computational complexity would be $O\left(N^{\text{fr}}_{\text{itr}}N_{\text{el}}N_pN_xN_yN_z\right)$ and the required memory would be proportional to $O\left(N_{\text{el}}N_pN_xN_yN_z\right)$. The obtained solution is the full-rank solution. 

\section{Proper Generalised Decomposition Solution}\label{Sec:PGD}

This section details the proper generalised decomposition strategies for \textit{a priori} low-rank solutions of the Boltzmann equation. To simplify the presentation, the subindex $a=\{P,T\}$ will be omitted for the perturbed velocity distribution function $h_a$.  

\subsection{Separation of spatial and velocity coordinates}

To construct a low-rank solution to problem \eqref{Shakhov} via PGD, the molecular velocity is expressed in the cylindrical coordinates as
\begin{equation}\label{cylind}
    v_x=v_r\cos\theta,\quad v_y=v_r\sin\theta,\quad v_z=v_z,
\end{equation}
where $\theta\in\left[0,2\pi\right]$ is the polar angle in the $x$-$y$ plane and $v_r=\sqrt{v^2_x+v^2_y}$. Substituting Eq.~\eqref{cylind} into Eqs.~\eqref{Shakhov} and~\eqref{moments}, the kinetic equation turns into
\begin{equation}\label{ShakhovCylinder}
    \begin{aligned}
        v_r\cos\theta\frac{\partial h}{\partial x}+v_r\sin\theta\frac{\partial h}{\partial y}=\delta\left[2uv_zf_{\text{eq}}+\frac{4}{15}qv_z\left(\bm{v}^2-\frac{5}{2}\right)f_{\text{eq}}-h\right]-s,\\
        \left(x,y\right)\in\Omega,\\
    h=0,\quad \left(x,y\right)\in\partial\Omega,\ (\cos\theta,\sin\theta)\cdot\bm{n}<0.
    \end{aligned}
\end{equation}
The independent variables of the system become $x$, $y$, $v_r$, $\theta$ and $v_z$.

\subsubsection{Separated representation of perturbed velocity distribution function}

PGD assumes that a low-rank solution is expressible in terms of a finite sum of products of functions that depend on parts of the coordinates, proposed as
\begin{equation}\label{PGDh}
    h\left(x,y,v_r,\theta,v_z\right)=\sum^m_{i=1}X_i\left(x,y,\theta\right)V_i\left(v_r,v_z\right),
\end{equation}
where $X_i$ and $V_i$ are unknown functions, named the $i$-th spatio-angular mode and the $i$-th mesoscopic velocity mode, respectively; and $m$ is the total number of PGD modes. The flow velocity and heat flux are, correspondingly, written in separated forms as
\begin{equation}
    u\left(x,y\right)=\sum^m_{i=1}Y_i\left(x,y\right)U_i,\quad q\left(x,y\right)=\sum^m_{i=1}Y_i\left(x,y\right)Q_i,
\end{equation}
where
\begin{equation}\label{PGDmoments}
    Y_i=\int X_i\mathrm{d}\theta,\ U_i=\int v_zV_iv_r\mathrm{d}v_r\mathrm{d}v_z,\ Q_i=\int v_z\left(\bm{v}^2-\frac{5}{2}\right)V_iv_r\mathrm{d}v_r\mathrm{d}v_z,
\end{equation}
are the $i$-th macroscopic spatial, velocity and heat flux modes. $Y_i$, $U_i$, $Q_i$ and their products are also named modes. In particular, $Y_i\times U_i$ and $Y_i\times Q_i$ are modes of the velocity and heat flux fields, respectively. Note that when constructing the separated representation~\eqref{PGDh}, the coordinates $x$, $y$ and $\theta$ are grouped. This is essential because the polar angle $\theta$ dominates the inflow/outflow at the boundary of a general two-dimensional spatial domain, where an appropriate boundary condition should be imposed to solve a boundary value problem with a hyperbolic operator. 


The PGD modes are commonly found by the Galerkin method and are enriched sequentially~\cite{PGDprimer}; that is, when computing the $m$-th modes, the first $m-1$ modes are known. Thus, the $m$-th modes $X_m$ and $V_m$ are found such that
\begin{equation}
    \begin{aligned}
        \int X^{\ast}V^{\ast}\left\{\delta X_mV_m+v_r\cos\theta V_m\frac{\partial X_m}{\partial x}+v_r\sin\theta V_m\frac{\partial X_m}{\partial y}\right.\\
        -2\delta Y_mU_mv_zf_{\text{eq}}-\frac{4}{15}\delta Y_mQ_mv_z\left(\bm{v}^2-\frac{5}{2}\right)f_{\text{eq}}\\
        = -s-\sum^{m-1}_{i=1}\left[\delta X_iV_i+v_r\cos\theta V_i\frac{\partial X_i}{\partial x}+v_r\sin\theta V_i\frac{\partial X_i}{\partial y}\right.\\
        \left.\left.-2\delta Y_iU_iv_zf_{\text{eq}}-\frac{4}{15}\delta Y_iQ_iv_z\left(\bm{v}^2-\frac{5}{2}\right)f_{\text{eq}}\right]\right\}\mathrm{d}x\mathrm{d}y\mathrm{d}\bm{v},
    \end{aligned}
\end{equation}
where $X^{\ast}V^{\ast}$ are some test functions. Although the original physical equation is linear, PGD results in a high-dimensional nonlinear problem. Fix-point iteration with an alternating-direction strategy is commonly used to approximate the nonlinear problem~\cite {PGDprimer}. Due to the separated representation, the computational cost is reduced since the high-dimensional problem is decomposed into a series of low-dimensional ones. Specifically, at each fixed-point iteration, the modes $X_m$ and $V_m$ are calculated one by one.  

\subsubsection{Mesoscopic velocity mode}

The goal is to compute the $m$-th velocity mode $V_m$ by assuming that $X_m$, as well as $Y_m$, are known. Thus, the test function is set to $X^{\ast}=X_m$. $V_m$ can be found from the following strong form
\begin{equation}\label{Vstrong}
    \begin{aligned}
        \left(\delta\alpha^m_m+v_r\beta^m_m \right)V_m-\delta\gamma^m_m\left[2U_mv_z-\frac{4}{15}Q_mv_z\left(\bm{v}^2-\frac{5}{2}\right)\right]f_{\text{eq}}=-\sigma^m s\\
         -\sum^{m-1}_{i=1}\left(\delta\alpha^m_i+v_r\beta^m_i\right)V_i-\delta\gamma^m_i\left[2U_iv_z-\frac{4}{15}Q_iv_z\left(\bm{v}^2-\frac{5}{2}\right)\right]f_{\text{eq}},
    \end{aligned}
\end{equation}
where
\begin{equation}
    \begin{aligned}
        \alpha^k_i=\int X_{k}X_{i}\mathrm{d}x\mathrm{d}y\mathrm{d}\theta,\quad
        \beta^k_{i}=\int\left(\cos\theta X_{k}\frac{\partial X_{i}}{\partial x}+\sin\theta X_{k}\frac{\partial X_{i}}{\partial y}\right) \mathrm{d}x\mathrm{d}y\mathrm{d}\theta,\\
        \gamma^k_i=\int X_{k}Y_{i}\mathrm{d}x\mathrm{d}y\mathrm{d}\theta,\quad\sigma^k=\int X_{k}\mathrm{d}x\mathrm{d}y\mathrm{d}\theta.
    \end{aligned}
\end{equation}

\subsubsection{Spatio-angular mode}

The goal is to compute the $m$-th spatio-angular mode $X_m$ after $V_m$, $U_m$, and $Q_m$ are computed. The test function is set to $V^{\ast}=V_m$. $X_m$ is then obtained from the following boundary value problem
\begin{equation}~\label{Xstrong}
    \begin{aligned}
        \delta\hat{\alpha}_m X_m+ \hat{\beta}_m\left(\cos\theta\frac{\partial X_m}{\partial x}+\sin\theta\frac{\partial X_m}{\partial y}\right)-2\delta\hat{\gamma}Y_mU_m-\frac{4}{15}\delta\hat{\kappa}Y_mQ_m \\
       = -\hat{\sigma}-\sum^{m-1}_{i=1}\left[\delta\hat{\alpha}_i X_i+\hat{\beta}_i\left(\cos\theta\frac{\partial X_i}{\partial x}+\sin\theta\frac{\partial X_i}{\partial y}\right)\right.\\
       \left.-2\delta\hat{\gamma}Y_iU_i-\frac{4}{15}\delta\hat{\kappa}Y_iQ_i\right],\quad\left(x,y\right)\in\Omega,\\
       X_m=0,\quad \left(x,y\right)\in\partial\Omega,\ (\cos\theta,\sin\theta)\cdot\bm{n}<0
    \end{aligned}
\end{equation}
where
\begin{equation}
    \begin{aligned}
        \hat{\alpha}_i=\int V_{m}V_{i}v_r\mathrm{d}v_r\mathrm{d}v_z,\quad \hat{\beta}_{i}=\int v_rV_{m}V_{i}v_r\mathrm{d}v_r\mathrm{d}v_z,\\
        \hat{\gamma}=\int v_zV_{m}f_{\text{eq}}v_r\mathrm{d}v_r\mathrm{d}v_z,\quad\hat{\kappa}=\int v_z\left(\bm{v}^2-\frac{5}{2}\right)V_{m}f_{\text{eq}}v_r\mathrm{d}v_r\mathrm{d}v_z,\\
        \hat{\sigma}=\int sV_{m}v_r\mathrm{d}v_r\mathrm{d}v_z.
    \end{aligned}
\end{equation}
Note that the problem described in \eqref{Xstrong} has a very similar structure to the kinetic equation~\eqref{ShakhovCylinder}, if $X_m$ and $Y_m$ are viewed as analogous to the velocity distribution function $h$ and the macroscopic quantities $u$ and $q$. Thus, its weak form is solved by the discontinuous Galerkin method.

\subsubsection{Stopping criteria and update of mesoscopic velocity modes}

A PGD algorithm is generally constructed as a two-loop structure~\cite{PGDprimer}: starting from the first mode, the outer loop enriches the PGD modes, and the inner loop applies alternative direction fixed-point iterations to search for the velocity and spatio-angular modes iteratively. To trade off the accuracy and computational cost, a stopping criterion should be given. In this work, the outer loop is restricted by a prescribed truncation number of modes $M_{\text{md}}$. The inner loop is stopped when the relative difference $\mathtt{res}$ of the amplitude of the PGD mode $A_m$ in two consecutive steps is less than a tolerance $\mathtt{tol}$, i.e.,
\begin{equation}
    \mathtt{res}=\frac{\vert A^{\mathtt{itr}}_m - A^{\mathtt{itr}-1}_m \vert}{A^{\mathtt{itr}-1}_m}<\mathtt{tol}
\end{equation}
or the iteration exceeds the allowed maximum step $N^{\text{in}}_{\text{itr}}$. Here, the superscript $^\mathtt{itr}$ is the index of the iteration step, and the amplitude is defined as
\begin{equation}
    A_m=\|X_m\|\cdot\|V_m\|,
\end{equation}
where $\|\cdot\|$ is the standard $\mathcal{L}_2$ norm.

As PGD enrichment progresses (i.e., $m$ increases), the convergence of the fix-point iteration may be difficult, especially when the rarefaction parameter $\delta$ is large, leading to a PGD solution with significant error. To overcome this defect, the velocity modes $\{V_i\}^m_{i=1}$ are updated when the fix-point iteration is completed. That is, $\{V_i\}^m_{i=1}$ are solved together from 
\begin{equation}\label{VstrongUpdate}
    \begin{aligned}
        \sum^{m}_{i=1}\left(\delta\alpha^k_i+v_r\beta^k_i\right)V_i-\delta\gamma^k_i\left[2U_iv_z-\frac{4}{15}Q_iv_z\left(\bm{v}^2-\frac{5}{2}\right)\right]f_{\text{eq}}=-\sigma^k s,\\
         k=1,\dots,m,
    \end{aligned}
\end{equation}
based on the obtained spatio-angular modes $\{X_i\}^m_{i=1}$. This PGD strategy with update was first proposed for time-dependent PDEs~\cite{NOUY20101603}, where numerical tests have demonstrated its effectiveness in improving convergence when determining PGD modes.  

\subsubsection{Algorithm and computational complexity}

The implementation of PGD is presented in \textbf{Algorithm}~\ref{Algorithm}. Note that the spatio-angular mode $X_m$ and the macroscopic spatial mode $Y_m$ are normalised in line 7 for numerical stability.

\begin{algorithm}[t]
\caption{PGD algorithm}\label{Algorithm}
\begin{algorithmic}[1]
\State Specify user-controlled input $M_{\text{md}}$, $N^{\text{in}}_{\text{itr}}$, $\mathtt{tol}$
\For{$m=1,\dots, M_{\text{md}}$}
    \State Initialise $X_m$ and $V_m$
    \State $k\leftarrow1$
    \While{$\mathtt{res}>\mathtt{tol}^{\text{in}}_{\text{itr}}$ and $k\leq N^{\text{in}}_{\text{itr}}$}
        \State Solve Eq.~\eqref{Vstrong} for $V_m$; Update $U_m$ and $Q_m$
        \State Solve Eq.~\eqref{Xstrong} for $X_m$; Update $Y_m$
        \State Normalised $X_m$ and $Y_m$: $X_m\leftarrow X_m/\| X_m\|$, $Y_m\leftarrow Y_m/\| Y_m\|$
        \State Update $A^{\mathtt{itr}}_m$; Calculate $\mathtt{res}$
        \State $k\leftarrow k+1$
    \EndWhile
    \State Solve Eq.~\eqref{VstrongUpdate} to update $\{V_i\}^m_{i=1}$
\EndFor
\end{algorithmic}
\end{algorithm}

When solving Eqs.~\eqref{Vstrong} and~\eqref{VstrongUpdate} for solutions of mesoscopic velocity modes, the discrete velocity method is still applicable. If $v_r$ and $v_z$ are discretised by $\{v^{j_1}_r\}^{N_r}_{j_1=1}$ and $\{v^{j_2}_z\}^{N_z}_{j_2=1}$ points, respectively, the macroscopic velocity and heat flux modes $U_i$ and $Q_i$ are approximated by quadrature 
\begin{equation*}
    U_i=\sum v^{j_2}_zV^{j}_iv^{j_1}_r\Delta v^{j},\quad Q_i=\sum v^{j_2}_z\left[\left(v^{j_1}_r\right)^2+\left(v^{j_2}_z\right)^2-\frac{5}{2}\right]V^{j}_iv^{j_1}_r\Delta v^{j},
\end{equation*}
where $V^{j}_i=V_i\left(v^{j_1}_r,v^{j_2}_z\right)$ and $\Delta v^{j}$ is the quadrature weights. The resultant linear system has a size of $N_r\times N_z$ or $m\times N_r\times N_z$, which can be handled by LU-decomposition. For spatio-angular modes, it is required to solve an equation [i.e., Eq.~\eqref{Xstrong}] with a structure very similar to the original kinetic equation. The same solution strategy used for the full-rank solution is adopted: the polar angle is discretised by $N_{\theta}$ points $\{\theta^{j_3}\}^{N_{\theta}}_{j_3=1}$, and the macroscopic spatial modes are estimated through the quadrature $Y_i=\sum X^{j_3}_i\Delta\theta^{j_3}$ with $X^{j_3}_i=X_i\left(x,y,\theta^{j_3}\right)$ and $\Delta\theta^{j_3}$ being the quadrature weight; $X^{j_3}_i$ is approximated using the discontinuous Galerkin method. $X^{j_3}_m$ is also solved iteratively. However, unlike for the full-rank solution, the number of iterative steps, denoted as $N^{\text{si}}_{\text{itr}}$, to obtain a convergent $X_m$ is almost independent of the value of the rarefaction parameter, and only a few steps (e.g., no greater than 9 for all tests presented in \S\ref{Sec:Result}) are needed to achieve convergence. 

In summary, the computational complexity of the PGD approximation is around $O\left(M_{\text{md}} N^{\text{in}}_{\text{itr}}\left(N_rN_z+N^{\text{si}}_{\text{itr}}N_{\text{el}}N_p N_{\theta}\right)+\left(M_{\text{md}}\right)!N_rN_z\right)$ and the memory requirement is proportional to $O\left(M_{\text{md}}\left(N_rN_z+N_{\text{el}}N_pN_{\theta}\right)\right)$. Compared to the full-rank solution, the computational complexity and the required memory can be significantly reduced as $M_{\text{md}}$ and $N^{\text{in}}_{\text{itr}}$ are expected to be small. This implies that a good approximation can be achieved with a few PGD modes, which are obtained from a quickly converged fixed-point iteration in the enrichment.

\subsection{Treating rarefaction parameter as an extra coordinate}

In practice, parametric analysis is typically performed by calculating flow fields over a wide range of rarefaction degrees. PGD provides a method for rapidly solving the parametrised kinetic equation, where the rarefaction parameter $\delta$ is considered as an additional coordinate. Flow properties at any point within the domain of the parameter can be accessed immediately, once the computational vademecum of flow velocity and heat flux is obtained.

The PGD solutions for $h$ now also depend on $\delta$, and can be expressed as
\begin{equation} h\left(x,y,v_r,\theta,v_z,\delta\right)=\sum^m_{i=1}X_i\left(x,y,\theta\right)V_i\left(v_r,v_z,\delta\right).
\end{equation}
The mesoscopic velocity modes are assumed to be parametrised with similar separated representations. The flow velocity and heat flux, thereby, read
\begin{equation}
    u\left(x,y,\delta\right)=\sum^m_{i=1}Y_i\left(x,y\right)U_i\left(\delta\right),\quad q\left(x,y,\delta\right)=\sum^m_{i=1}Y_i\left(x,y\right)Q_i\left(\delta\right),
\end{equation}
where the macroscopic modes are still calculated from Eq.~\eqref{PGDmoments}. In contrast to the widely used strategy of introducing parametric modes (additional modes as functions of $\delta$)~\cite{SEVILLA2020112631}, parametrising existing modes $V_i$, $U_i$, and $Q_i$ can keep the greedy PGD solution process as concise as possible. The algorithm for obtaining the PGD modes is the same as the Algorithm~\ref{Algorithm}, except for the following modifications. 

Rather than a single value input, the parameter $\delta$ in Eqs.~\eqref{Vstrong} and~\eqref{VstrongUpdate} is an independent variable that covers a continuous range. When solving the equations, it is discretised by $N_{\delta}$ discrete nodes $\{\delta^j\}^{N_{\delta}}_{j=1}$. Consequently, the sizes of the linear systems for finding $V_m$ and updating $\{V_i\}^m_{i=1}$ become $N_r\times N_z\times N_{\delta}$ and $m\times N_r\times N_z\times N_{\delta}$, respectively. Due to the dependency on the parameter, the boundary value problem for solutions of spatio-angular modes is modified as
\begin{equation}~\label{XstrongPara}
    \begin{aligned}
       \tilde{\alpha}_m X_m+ \tilde{\beta}_m\left(\cos\theta\frac{\partial X_m}{\partial x}+\sin\theta\frac{\partial X_m}{\partial y}\right)-2\tilde{\gamma}_mY_m-\frac{4}{15}\tilde{\kappa}_mY_m \\
       = -\tilde{\sigma}-\sum^{m-1}_{i=1}\left[\tilde{\alpha}_i X_i+\tilde{\beta}_i\left(\cos\theta\frac{\partial X_i}{\partial x}+\sin\theta\frac{\partial X_i}{\partial y}\right)\right.,\\
       \left.-2\tilde{\gamma}_iY_i-\frac{4}{15}\tilde{\kappa}_iY_i\right],\quad\left(x,y\right)\in\Omega,\\
       X_m=0,\quad \left(x,y\right)\in\partial\Omega,\ (\cos\theta,\sin\theta)\cdot\bm{n}<0
    \end{aligned}
\end{equation}
where
\begin{equation}
    \begin{aligned}
        \tilde{\alpha}_i=\int\delta V_{m}V_{i}v_r\mathrm{d}v_r\mathrm{d}v_z\mathrm{d}\delta,\quad \tilde{\beta}_{i}=\int v_rV_{m}V_{i}v_r\mathrm{d}v_r\mathrm{d}v_z\mathrm{d}\delta,\\
        \tilde{\gamma}_i=\int \delta v_zV_{m}U_if_{\text{eq}}v_r\mathrm{d}v_r\mathrm{d}v_z\mathrm{d}\delta,\\
        \tilde{\kappa}_i=\int\delta v_z\left(\bm{v}^2-\frac{5}{2}\right)V_{m}Q_if_{\text{eq}}v_r\mathrm{d}v_r\mathrm{d}v_z\mathrm{d}\delta,\quad\tilde{\sigma}=\int sV_{m}v_r\mathrm{d}v_r\mathrm{d}v_z\mathrm{d}\delta.
    \end{aligned}
\end{equation}
But the solution strategy remains unchanged. 

For the parametric problem, the computational complexity and required memory increase to $O\left(M_{\text{md}}N^{\text{in}}_{\text{itr}}\left(N_r N_zN_{\delta}+N^{\text{si}}_{\text{itr}}N_{\text{el}}N_p N_{\theta}\right)+\left(M_{\text{md}}\right)!N_rN_zN_{\delta}\right)$ and $O\left(M_{\text{md}}\left(N_rN_zN_{\delta}+N_{\text{el}}N_pN_{\theta}\right)\right)$, respectively.

\section{Numerical examples and Discussions}\label{Sec:Result}

\begin{figure}[t]
\centering
\includegraphics[width=\textwidth]{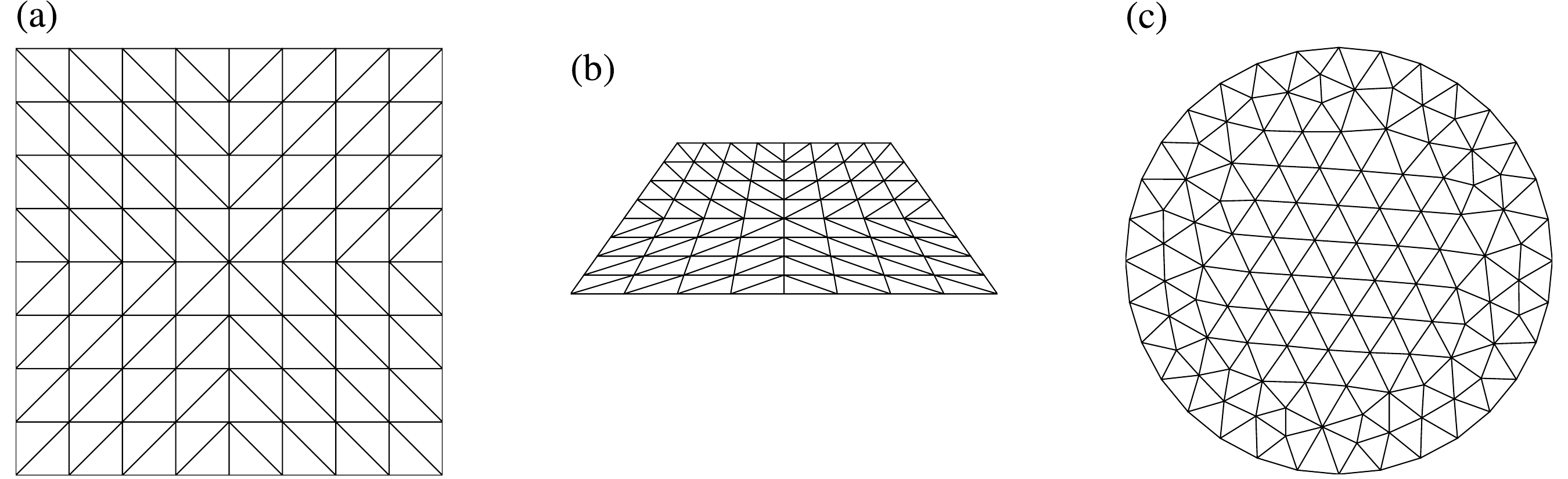}
\caption{Geometry and triangular mesh for two-dimensional spatial domain. (a) Squara domain partitioning by 128 triangular elements; (b) Isosceles trapezoidal domain with a ratio of bases of 0.5 and an acute angle of 54.74\textdegree, partitioning by 128 elements; (c) Circular domain partitioning by 780 elements.}\label{Mesh}
\end{figure}

In this section, numerical examples of Poiseuille and thermal creep flows through long straight channels with square, trapezoidal, and circular cross sections are presented to demonstrate the performance of PGD solutions for this high-dimensional and parametrised problem governed by kinetic equations. Figure~\ref{Mesh} shows the two-dimensional spatial domains and the triangular meshes for the discontinuous Galerkin finite element approximation. The side length, lower base, and radius of the three cross sections are set as the characteristic flow length. The meshes contain 128, 128, and 780 triangles, respectively. The degree of approximation polynomial basis selected for the spatial elements is $p=3$. The truncated velocity domains $v_r\in[0,4]$ and $v_z\in[-4,4]$ are discredited by $N_r=N_z=24$ non-uniformly distributed nodes~\cite{SUImplicitDG}. The polar angle $\theta\in[0,2\pi]$ is divided into $N_{\theta}=48$ intervals of equal spacing. Integrations involving velocity and angle are approximated by using the midpoint rule. For comparison, the same resolutions are used for the full-rank results. When considering the parametric problem, the range of the rarefaction parameter $\delta\in[0.01,100]$ is discretised by $N_{\delta}=33$ nodes, and the corresponding integration is evaluated by Simpson's rule. For PGD solutions, the maximum number of modes is set to $M_{\text{md}}=15$, and the maximum allowed iteration steps and the tolerance for the fixed-point iteration are $N^{\text{in}}_{\text{itr}}=10$ and $\mathtt{tol}=10^{-3}$. All tests are run on a single Intel\textregistered~core\texttrademark~i7-12700K CPU.

\subsection{Flow through channel with square cross section}

\begin{figure}[t]
\centering
\includegraphics[width=\textwidth]{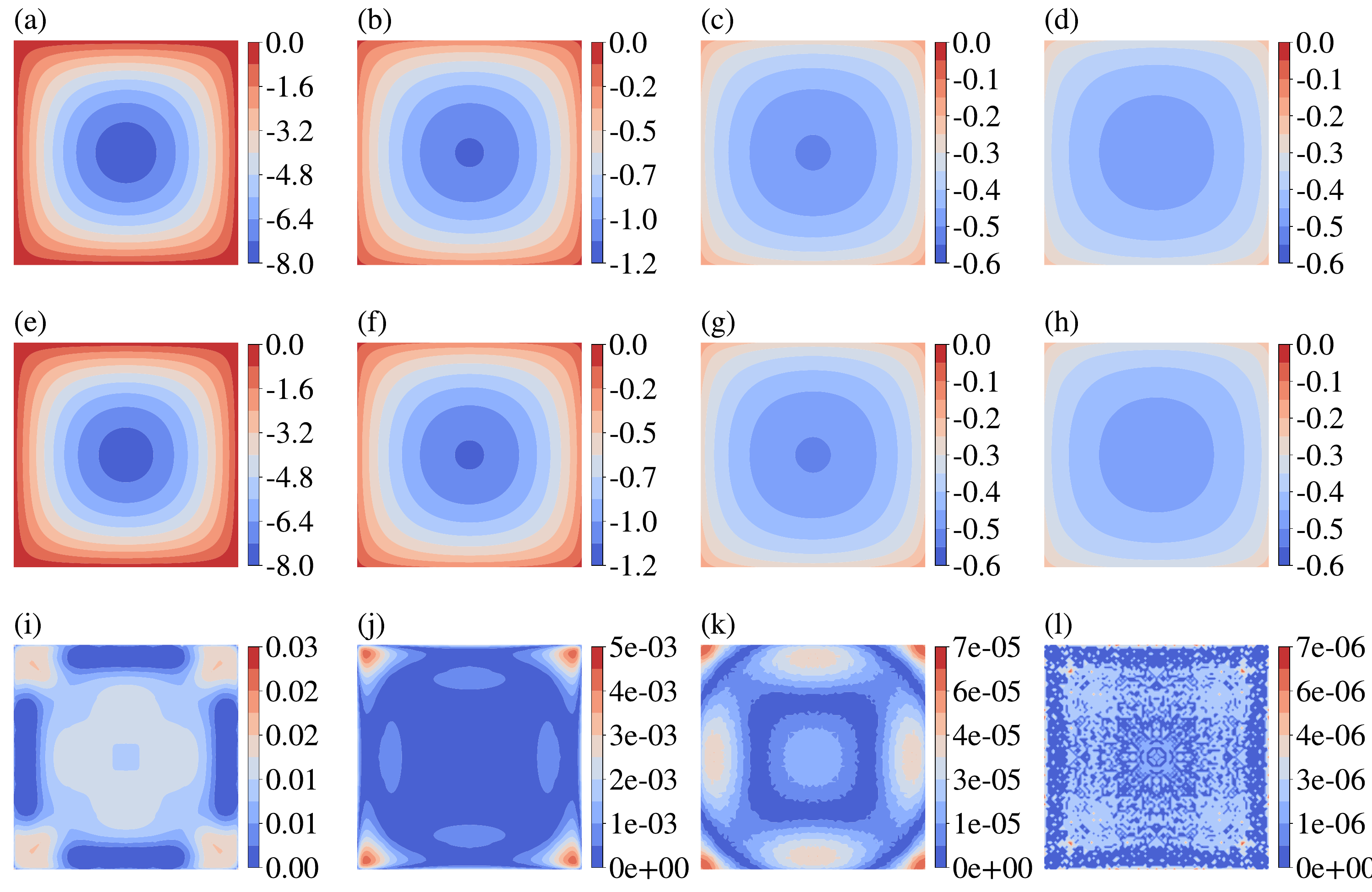}
\caption{Flow velocity of Poiseuille flow through a long channel with square cross section. First row (a)-(d): PGD solutions with 15 modes; second row (e)-(h): full-rank solutions; third row (i)-(l): relative errors between the low- and full-rank solutions. First column (a), (e) and (i): $\delta=100$; second column (b), (f) and (j): $\delta=10$; third column (c), (g) and (k): $\delta=1$; fourth column (d), (h) and (l): $\delta=0.1$.}\label{Ufield}
\end{figure}

\begin{figure}[h]
\centering
\includegraphics[width=\textwidth]{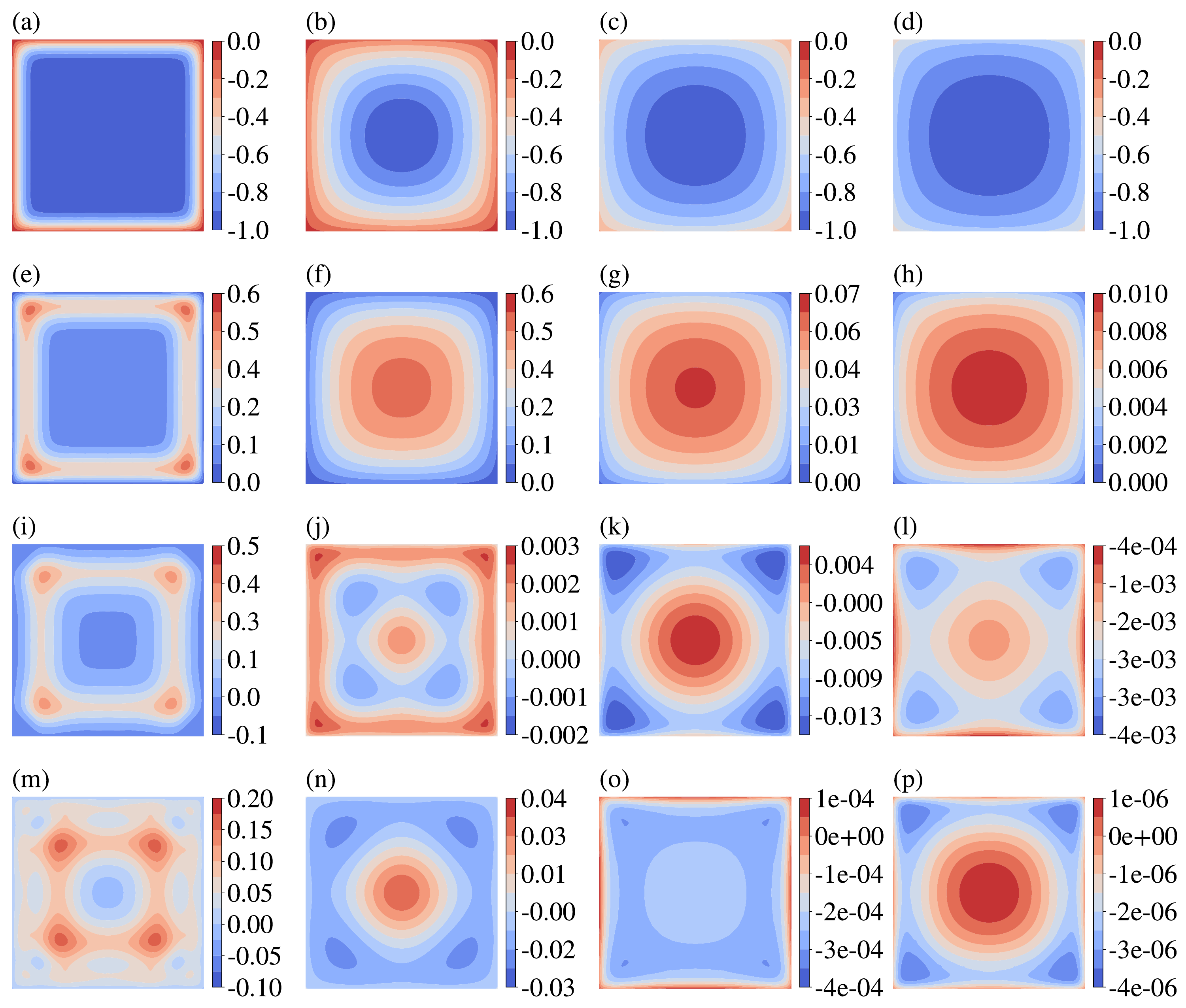}
\caption{First four normalised PGD modes of velocity field $Y_i\times U_i/\max\{|Y_1\times U_1|\}$ of Poiseuille flow through a long channel with square cross section. First row (a)-(d): $i=1$; second row (e)-(h): $i=2$; third row (i)-(l): $i=3$; fourth row (m)-(p): $i=4$. First column (a), (e), (i) and (m): $\delta=100$; second column (b), (f), (j) and (n): $\delta=10$; third column (c), (g), (k) and (o): $\delta=1$; fourth column (d), (h), (l) and (p): $\delta=0.1$.}\label{Umode}
\end{figure}

Figure~\ref{Ufield} illustrates the macroscopic flow velocity of the Poiseuille flow through a long channel with a square cross section. Comparison between PGD solutions and full-rank solutions is visualised with the contour of relative error $|u_{\text{PGD}}-u_{\text{full-rank}}|/|u_{\text{full-rank}}|$. Four different rarefactions $\delta=100$, $\delta=10$, $\delta =1$ and $\delta=0.1$ are considered. It shows that PGD with 15 modes recovers the flow field accurately to within $3\%$. In particular, for smaller $\delta$, i.e., when the flow is more rarefied, PGD solution demonstrates its effectiveness with a very high accuracy. Figure~\ref{Umode} plots the first four modes of the macroscopic velocity field normalised by the maximum magnitude of the first mode, namely $Y_i\times U_i/\max\{|Y_1\times U_1|\}$. It can be seen that the first mode captures the most relevant and global feature of the solution, while the next modes introduce local corrections, just as usual in reduced-order models~\cite{SEVILLA2020112631}.

\begin{figure}[h]
\centering
\includegraphics[width=\textwidth]{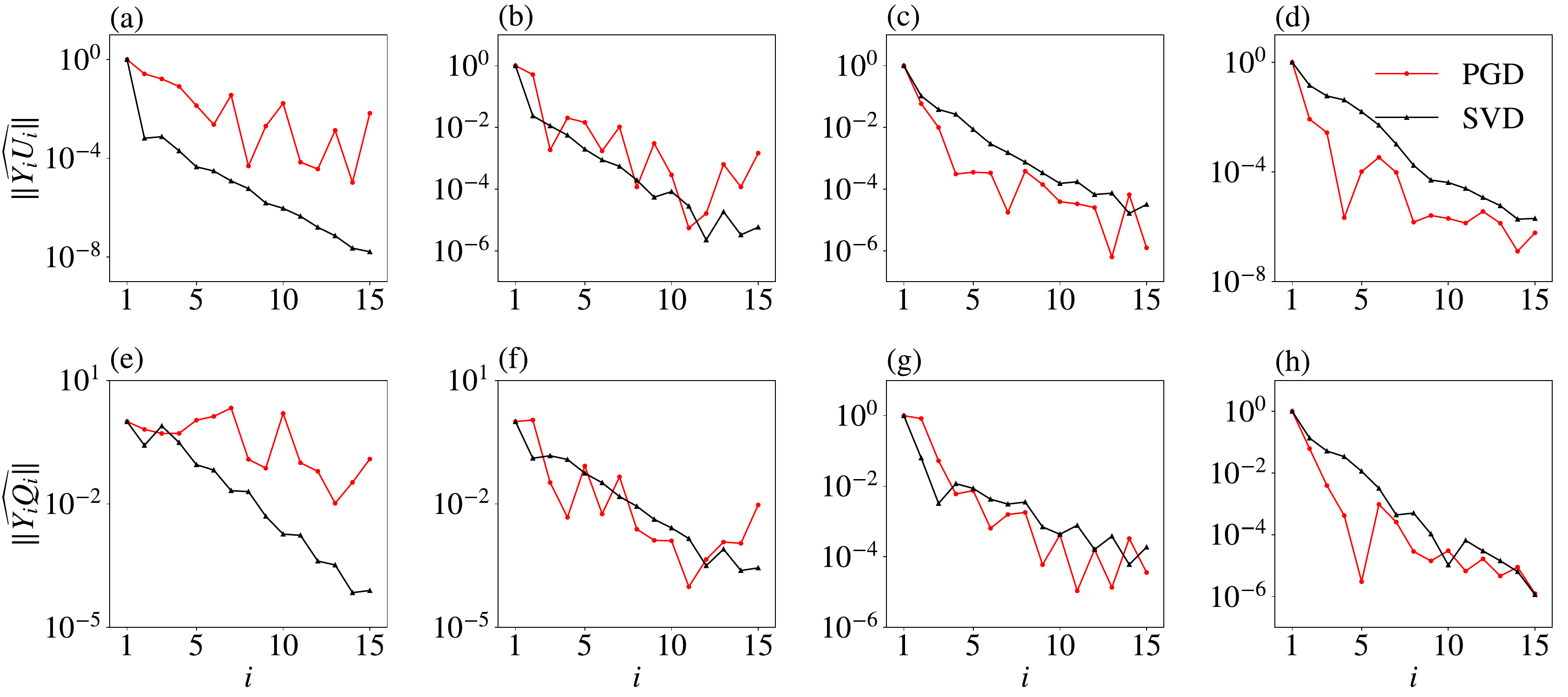}
\caption{Relative amplitudes of the modes of velocity and heat flux field: $\|\widehat{Y_iU_i}\|=\|Y_i\|\cdot\|U_i\|/\|Y_1\|\cdot\|U_1\|$ and $\|\widehat{Y_iQ_i}\|=\|Y_i\|\cdot\|Q_i\|/\|Y_1\|\cdot\|Q_1\|$. Red lines with circular markers show PGD modes. Black lines with triangular marks are the modes obtained by SVD for the full-rank solutions. First column (a) and (e): $\delta=100$; second column (b) and (f): $\delta=10$; third column (c) and (g): $\delta=1$; fourth column (d) and (h): $\delta=0.1$. }\label{modeAmp}
\end{figure}

To quantify the importance of the modes in the PGD solution, the relative amplitudes (related to those of the first) of the modes of the velocity and heat flux fields are plotted in Figure~\ref{modeAmp}. The results are compared with the modes obtained by applying SVD to the full-rank solution, where SVD provides the optimal choice for the reduced basis~\cite{PGDprimer}. It is observed that the SVD results decay dramatically versus the modes, indicating the low-rank structure of the solutions. For $\delta=10$, 1 and 0.1, the magnitude and downward trend of the PGD mode amplitudes are very close to those of the SVD results, showing a good performance of the PGD algorithm in searching for the reduced basis. When $\delta=100$, the amplitudes of the PGD modes generally drop, although they decay slowly and can undergo large oscillations. Note that when the gas is slightly rarefied, the right-hand side of the kinetic equation becomes stiff. This may cause the fixed-point iteration having difficulty to converge, leading to a relatively large error in the PGD solutions with a few modes; see Figure~\ref{Ufield}(i).           

\begin{table}[t]
\centering
\begin{tabular}{ccccccccc}
\hline
\multirow{2}{*}{$\delta$} & \multicolumn{4}{c}{Full-rank} & & \multicolumn{3}{c}{PGD} \\
\cline{2-5}\cline{7-9}
 & $G_P$ & $G_T$ & $N^{\text{fr}}_{\text{itr}}$ & $t_{\text{CPU}}$, [h] & & $G_P$ & $G_T$ & $t_{\text{CPU}}$, [h] \\
\hline
0.1 & 0.7950 & 0.3648 &  6 & 0.007 & & 0.7950 & 0.3648 & 0.008\\
0.5 & 0.7675 & 0.2963 & 11 & 0.014 & & 0.7675 & 0.2963 & 0.009\\
1.0 & 0.7755 & 0.2557 & 16 & 0.020 & & 0.7755 & 0.2557 & 0.009\\
2.0 & 0.8195 & 0.2083 & 25 & 0.031 & & 0.8189 & 0.2083 & 0.009\\
5.0 & 0.9975 & 0.1379 & 61 & 0.076 & & 0.9976 & 0.1379 & 0.010\\
10.0& 1.3258 & 0.0879 & 144& 0.174 & & 1.3258 & 0.0880 & 0.010\\
20.0& 2.0086 & 0.0503 & 391& 0.466 & & 2.0095 & 0.0504 & 0.010\\
40.0& 3.3952 & 0.0271 & 1146& 1.349& & 3.4087 & 0.0266 & 0.010\\
50.0& 4.0891 & 0.0220 & 1630& 1.918& & 4.1020 & 0.0218 & 0.010\\
80.0& 6.1589 & 0.0140 & 3411& 4.007& & 6.1848 & 0.0145 & 0.010\\
100.0& 7.5221& 0.0113 & 4811& 5.639& & 7.5997 & 0.0114 & 0.010\\
\hline
\end{tabular}
\caption{Comparison of full-rank and PGD solutions at different rarefaction parameters: $G_P$ -- Poiseuille flow rate; $G_T$ -- thermal creep flow rate; $N^{\text{fr}}_{\text{itr}}$ -- number of iterations to obtain converged full-rank solutions (maximum relative difference in flow velocity and heat flux between two successive iteration steps is less than $10^{-5}$); $t_{\text{CPU}}$ -- CPU time measured in hour.}\label{Table1}
\end{table}

To show the efficiency in finding low-rank PGD solutions, the CPU times to obtain full-rank and PGD solutions are compared for different rarefaction parameters in Table~\ref{Table1}. The flow rates of the Poiseuille and thermal creep flows are also listed for comparison. The PGD can produce almost the same flow rates as full-rank results when the flow is rarefied. The precision of the low-rank solutions slightly deteriorated when $\delta$ became large; however, the maximum difference (appearing in $G_T$ at $\delta=80$) is still less than 3.6\%. The CPU time to obtain a full-rank solution increases dramatically as $\delta$ increases, since the number of iterations to achieve convergence increases significantly~\cite{SU2020109245}. In contrast, PGD costs almost the same time to obtain a solution at any $\delta$, that is, the computational complexity of PGD remains the same for the whole range of rarefaction. Due to the reduced complexity, PGD can be much faster than the full-rank method, and the memory consumption by PGD is only about 2.6\% of that of the full-rank method.    

\clearpage
\subsection{Flow through channel with trapezoidal/circular cross section}

\begin{figure}[h]
\centering
\includegraphics[width=0.85\textwidth]{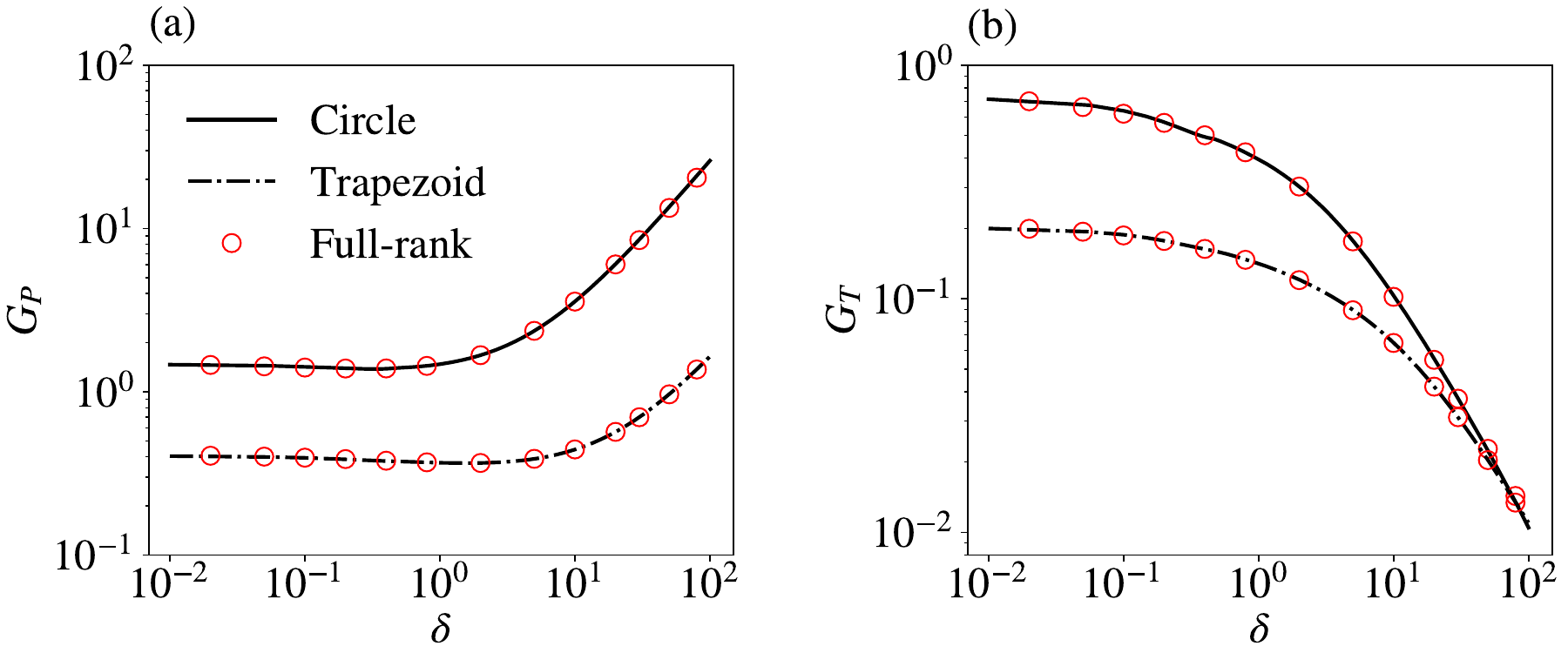}
\caption{(a) Dimensionless flow rate of Poiseuille flow and (b) dimensionless flow rate of thermal creep flow through long channels with trapezoidal and circular cross sections. Lines are PGD general solutions as continuous functions of rarefaction parameter $\delta$. Markers are full-rank solutions at some values of $\delta$.}\label{FRs}
\end{figure}

The parametrised kinetic equation~\eqref{ShakhovCylinder} with $\delta$ serving as an additional coordinate is solved by PGD for flows through long channels with trapezoidal and circular cross sections. Running the PGD calculation once, a general solution of the flow field is obtained that contains all possible data for $\delta$ varying from 0.01 to 100. Figure~\ref{FRs} illustrates dimensionless flow rates as continuous functions of the rarefaction parameter $\delta$. Full-rank solutions calculated at some values of $\delta$ are also included for comparison. PGD is found to be able to recover solutions with high precision in the whole range of $\delta$ considered. In addition to this, PGD only used 0.261 hours to obtain the general solution for the channel with the trapezoidal cross section and 0.405 hours for the channel with the circular cross section.     

Once the PGD general solution is obtained and stored as a computational vacuum, it can be used for many purposes, such as parametric analysis and the investigation of inverse problems. Here, as an example, the thermomolecular pressure difference (TPD) problem is considered. It is a specific state where the Poiseuille and thermal creep flows counterbalance and the net flow rate vanishes. The pressures and temperatures in the two reservoirs connected by the long channel satisfy the following relationship
\begin{equation}
    \frac{P_2}{P_1}=\left(\frac{T_2}{T_1}\right)^{\eta},
\end{equation}
where $\eta$ is the TPD coefficient. To evaluate this coefficient, it is assumed that the temperature ratio $T_2/T_1$ and the rarefaction $\delta_1$ with respect to $P_1$ and $T_1$ are known. The pressure ratio $P_2/P_1$ is solved from~\cite{Graur2009,Ritos01112011}
\begin{equation}\label{TPD}
    \frac{\mathrm{d}\tilde{P}}{\mathrm{d}\tilde{T}}=\frac{\tilde{P}}{\tilde{T}}\frac{G_T\left(\delta\right)}{G_P\left(\delta\right)},
\end{equation}
where $\tilde{P}=P/P_1$ and $\tilde{T}=T/T_1$ are the normalised local pressure and temperature along the flow direction, and $\delta=\delta_1\tilde{P}/\tilde{T}$ is the local rarefaction parameter. Using the PGD solutions of $G_P\left(\delta\right)$ and $G_T\left(\delta\right)$, $P_2/P_1$ can be obtained by integrating Eq.~\eqref{TPD} along $1\le\tilde{T}\le T_2/T_1$. The TPD coefficients for $T_2/T_1=3.8$ and various $\delta_1$ are listed in Table~\ref{Table2}. The results obtained from the general PGD solutions are compared with those in previous works based on full-rank calculations. It further shows the accuracy of the parametrised solutions obtained by PGD.

\begin{table}[t]
\centering
\begin{tabular}{cccccc}
\hline
$\delta_1$ & Circle & Circle in~\cite{Graur2009} & $\delta^{\ast}_1$ & Trapezoid & Trapezoid in~\cite{Ritos01112011}\\
\hline
0.02 & 0.4862 & 0.4849 & 0.02 & 0.4912 & 0.4907 \\
0.05 & 0.4765 & 0.4706 & 0.05 & 0.4815 & 0.4791 \\
0.1  & 0.4611 & 0.4531 & 0.1  & 0.4670 & 0.4646 \\
0.2  & 0.4330 & 0.4275 & 0.2  & 0.4441 & 0.4441 \\
0.5  & 0.3727 & 0.3761 & 0.5  & 0.4033 & 0.4019 \\
0.8  & 0.3366 & 0.3389 & 0.8  & 0.3725 & 0.3722 \\
1.0  & 0.3174 & 0.3186 & 1.0  & 0.3572 & 0.3560 \\
2.0  & 0.2458 & 0.2468 & 2.0  & 0.2970 & 0.2961 \\
5.0  & 0.1409 & 0.1423 & 5.0  & 0.2000 & 0.1991 \\
10.0 & 0.0739 & 0.0743 & 10.0 & 0.1253 & 0.1246 \\
20.0 & 0.0310 & 0.0306 & 20.0 & 0.0642 & 0.0647 \\
\hline
\end{tabular}
\caption{TPD coefficients for $T_2/T_1=3.8$ and various $\delta_1$. Results from PGD parametrised solutions are compared with those in previous works based on full-rank calculations. $\delta^{\ast}_1=0.4483\delta_1$ is the rarefaction parameter based on the hydraulic diameter of the channel.}\label{Table2}
\end{table}

\clearpage

\section{Conclusions}\label{Sec:conclusion}

A reduced-order solution strategy based on the PGD method is proposed for the high-dimensional parameterised Shakhov model equation, which describes the dynamics of rarefied gas flows driven by pressure and temperature gradients through a long channel with an arbitrary cross section. The low-rank solution of the velocity distribution function is approximated as the sum of a small number of function products, where each product contains a spatio-angular mode and a mesoscopic velocity mode. A sophisticated algorithm is designed to find these PGD modes. Due to the separated representation, evaluating the velocity distribution function is transformed into a few low-dimensional problems. Consequently, computational complexity is significantly reduced. The numerical results show that the method can obtain a solution in 1 minute at any Knudsen number and cost only 2.6\% of memory for the full-rank solution. For parametrised flows, the rarefaction parameter serves as an additional coordinate. A general solution is calculated once for all in the whole parameter range, and a specific solution is accessible in real time. Compared to full-rank solutions, the PGD solution possesses high accuracy.

The proposed reduced-order method is applicable for modelling particle transport processes beyond gas molecules, e.g., phonons, neutrons, photons, which are described by a kinetic equation.



\bibliographystyle{elsarticle-num} 
\bibliography{Ref}

\end{document}